\documentclass{amsart}
\usepackage{amssymb}
\usepackage{mathtools}
\usepackage{tikzit}
\tikzstyle{esquare}=[fill=white, draw=black, shape=square]
\tikzstyle{inny}=[shape=circle, fill=black, inner sep=0pt, minimum size=5pt]
\tikzstyle{ecircle}=[fill=white, draw=black, shape=circle]
\tikzstyle{znajdz}=[fill=none, draw=none, shape=square]
\tikzstyle{redy}=[fill=red, draw=none, shape=circle, inner sep=0pt, minimum size=5pt]

\tikzstyle{standard}=[-, draw=black, fill=none]
\tikzstyle{dashe}=[-, fill=none, dashed]
\tikzstyle{leftarow}=[<-]
\tikzstyle{righarow}=[->]
\tikzstyle{rede}=[-, draw=red, fill=none]

\usepackage{tikz}
\usepackage{tikz-cd}
\usetikzlibrary{arrows.meta}

\definecolor{section_col}{RGB}{0, 79, 159}
\definecolor{cite_col}{RGB}{0, 63, 127}

\usepackage{hyperref}
\hypersetup{
    colorlinks=true,
    linkcolor=section_col,
    filecolor=magenta,      
    urlcolor=cyan,
    citecolor=cite_col
    }
    
\newtheorem{Thm}{Theorem}[section]

\newtheorem{Cor}[Thm]{Corollary}
\newtheorem{Obs}[Thm]{Observation}
\newtheorem{Lem}[Thm]{Lemma}
\newtheorem{Fact}[Thm]{Fact}

\newtheorem{Prop}[Thm]{Proposition}
\theoremstyle{definition}%boldface title, romand body
\newtheorem{Rem}[Thm]{Remark}      
\theoremstyle{definition}
\newtheorem{Defn}[Thm]{Definition}
\newtheorem{Ex}[Thm]{Example}
\theoremstyle{definition}%boldface title, romand body
      
\theoremstyle{definition}
      
\theoremstyle{definition}
\newtheorem*{Conv}{Convention}
\theoremstyle{definition}

\newcommand{\thm}[1]{\begin{Thm} #1 \end{Thm}}
\newcommand{\lem}[1]{\begin{Lem} #1 \end{Lem}}

\newcommand{\prop}[1]{\begin{Prop} #1 \end{Prop}}

\newcommand{\defn}[1]{\begin{Defn} #1 \end{Defn}}

\newcommand{\pr}[1]{\begin{proof} #1 \end{proof}}

\newcommand{\R}{\mathbb{R}}  
\newcommand{\Z}{\mathbb{Z}}  
 
\newcommand{\C}{\mathbb{C}} 
\newcommand{\Ha}{\mathbb{H}}

\newcommand{\RP}{\R \mathrm{P}}
\newcommand{\CP}{\C \mathrm{P}}

\newcommand{\Spin}[1]{\text{Spin}(#1)}
\newcommand{\Sp}[1]{\text{Sp}(#1)}
\newcommand{\SU}[1]{\text{SU}(#1)}
\newcommand{\U}[1]{\text{U}(#1)}
\newcommand{\SO}[1]{\text{SO}(#1)}

\newcommand{\GL}[1]{\text{GL}(#1)}
\newcommand{\BSpin}[1]{B\text{Spin}(#1)}
\newcommand{\BSp}[1]{B\text{Sp}(#1)}

\newcommand{\ESpin}[1]{E\text{Spin}(#1)}
\newcommand{\BSU}[1]{B\text{SU}(#1)}
\newcommand{\BU}[1]{B\text{U}(#1)}
\newcommand{\BSO}[1]{B\text{SO}(#1)}

\newcommand{\ov}[1]{\overline{#1}}

\newcommand{\tylda}[1]{\widetilde{#1}}
\newcommand{\rksp}{\tylda{KSP}}
\newcommand{\rko}{\tylda{KO}}
\newcommand{\rku}{\tylda{K}}

\newcommand{\pso}[1]{P_{\text{SO}}( #1 )}
\newcommand{\pspin}[1]{P_{\text{Spin}}( #1 )}
\newcommand{\Cl}{\text{Cl}}

\title{Real spin bundles over $\CP^3$ and a new Euclidean embedding of $\RP^7$}
\author{Dominik Gdesz}
\date{\today}

\begin{document}
\begin{abstract}
    We generalize the $\alpha$-invariant introduced by Atiyah and Rees to an invariant of real spin bundles and use it to classify real bundles over $\CP^3$ admitting spin structure. We apply this result to show that $\RP^7$ can be smoothly embedded in $\R^{11}$.
\end{abstract}
\maketitle
\vspace{1cm}
\section{Introduction}
    Atiyah and Rees proved in \cite{Atiyah} that every continuous $\U2$ - bundle over $\CP^3$ admits a holomorphic structure. Their idea was to classify all such bundles and compare this classification with examples of holomorphic bundles that were already known at that time. They divided the problem of such classification into two cases depending on the parity of the first Chern class $c_1$. The case of bundles with even $c_1$ can be further reduced to $\SU2$ - bundles by tensoring with a tautological bundle. Atiyah and Rees used the fact that $\SU2\cong \Sp1$ to classify $\Sp1$ bundles instead, which can be achieved by passing to $KSP$ - theory. \\
    \indent Since $\Spin 3\cong \SU 2$, this result also classifies $\Spin 3$-bundles over $\CP^3$. We modify this approach in Section \ref{sec:inv} and introduce a general way of constructing invariants of spin bundles by using representations of spin groups coming from their actions on Clifford algebras and taking $K$-theoretic classes of bundles associated with these representations. These invariants completely classify $\Spin n$-bundles over connected CW-complexes of dimension $k$, provided that $n,k\leq 6$. On the other hand, bundles of dimension $\geq7$ over a CW-complex of dimension $\leq 6$ are stable. We use these facts to obtain the complete classification of real spin bundles over $\CP^3$, which is given in Section \ref{sec:main}. \\
    \indent In Section \ref{sec:app} we use this classification to prove that $\RP^7$ can be smoothly into $\R^{11}$. Since Hantzsche's argument in \cite{Hantzsche} it is known that $\RP^{7}$ cannot be embedded in $\R^{8}$. On the other hand, Mahowald proved in \cite{Mahowald} that $\RP^{7}$ embeds topologically in $\R^{11}$ and used a result by Haefliger from \cite{Haefliger1961/62} to conclude that it embeds smoothly in $\R^{12}$. This paper is the first improvement in estabilishing maximal codimension of a smooth Euclidean embedding for $\RP^7$ since Mahowald's result.

\newpage

\section{\label{sec:inv} Invariants of real spin bundles}

Recall that the group $\Spin n$ is contained in the even part $\Cl_n^0$ of the Clifford algebra $\Cl_n$ and that $\Cl_n^0\cong \Cl_{n-1}$. The irreducible representations of the algebras $\Cl_n$ were classified in the classical work of Atiyah, Bott, and Shapiro \cite{ATIYAH19643}. We summarize it briefly: 

\prop{There is one irreducible real representation of $\Cl_n$ if $n\not\equiv 3 \mod 4 $ and two real irreducible representations of equal dimensions if $n\equiv 3 \mod 4$. Moreover, for $n\equiv 1,5\mod 8$ these representations are automatically complex and for $n\equiv 2,3,4 \mod 8$ these representations are automatically quaternionic.}

We can restrict these representations to $\Spin n\subset \Cl_{n-1}$ and obtain representations of spin groups, which will be denoted as $\rho_n$ (or $\rho_n^1, \ \rho_n^2$ if $n=4k$). We will write $V_n$ for the real, complex, or quaternionic vector space on which $\Spin n$ acts via $\rho_n$ or $\rho_n^i$. Now suppose that we are given a principal $\Spin n$-bundle $P$ over some base space $X$. We can construct the associated bundle $P\times_{\rho_{n}}V_{n}$ (or two bundles $P\times_{\rho_{n}^i}V_{n}$ for $i=1,2$ if $n=4k$). If $n=2,6\mod 8$ then this bundle is equipped with a complex structure and if $n=3,4,5\mod 8$ then it is equipped with a quaternionic structure. We can now take an appropriate $K$-theoretical class and obtain an invariant of our spin bundle.
\defn{Let $P$ be a $\Spin n$ bundle over base $X$. We define the invariant $\rho[P]$ as:
\begin{align*}
[P\times_{\rho_{n}}V_{n}]&\in \rko(X) &\text{ if }\ n=8k+1,\\
[P\times_{\rho_{n}}V_{n}]&\in \rku(X) &\text{ if }\ n=8k+2, \\
[P\times_{\rho_{n}}V_{n}]&\in \rksp(X) &\text{ if }\ n=8k+3, \\
[P\times_{\rho_{n}^1}V_{n}]\oplus[P\times_{\rho_{n}^2}V_{n}]&\in \rksp(X)\oplus \rksp(X) &\text{ if }\ n=8k+4, \\
[P\times_{\rho_{n}}V_{n}]&\in \rksp(X) &\text{ if }\ n=8k+5, \\
[P\times_{\rho_{n}}V_{n}]&\in \rku(X) &\text{ if }\ n=8k+6, \\
[P\times_{\rho_{n}}V_{n}]&\in \rko(X) &\text{ if }\ n=8k+7, \\
[P\times_{\rho_{n}^1}V_{n}]\oplus [P\times_{\rho_{n}^2}V_{n}]&\in \rko(X)\oplus \rko(X) &\text{ if }\ n=8k+8. \\
\end{align*}
}
\lem{\label{lem:classification} If $X$ is a connected CW-complex of dimension $\leq 6$, then the invariant $\rho[P]$ defines a correspondence between $\Spin n$-bundles over $X$ for $n=3,4,5$ and appropriate $K$-theoretical groups. When $n=6$ it defines a correspondence between $\Spin 6$-bundles over $X$ and the image of $[X,\text{BSU}]$ in $\rku(X)$, i.e. the subgroup of $\rku(X)$ generated by bundles with $c_1=0$.}
\pr{The representations $\rho_n$ for $n =3,5,6$ and $\rho_4^i$ for $i=1,2$ give us the exceptional isomorphisms:
\begin{align*}
    \Spin 3 &\cong  \Sp1, \\
    \Spin 4 &\cong \Sp1\times\Sp1, \\
    \Spin 5 &\cong \Sp2, \\
    \Spin 6 &\cong \SU 4.
\end{align*}
Looking at the fibrations $\Sp n\hookrightarrow \Sp{n+1}\to S^{4n+3}$ and $\SU n\hookrightarrow \SU{n+1}\to S^{2n+1}$ we see that the inclusions $\BSp 1\hookrightarrow \BSp2\hookrightarrow \text{BSp}$ and $\BSU{4}\hookrightarrow \text{BSU}$ are $7$-equivalences, i.e. they induce isomorphisms on $\pi_i:i< 7$ and surjections on $\pi_7$. It follows that, since $X$ is a CW-complex of dimension $\leq 6$, there are set isomorphisms:
\begin{align*}
    [X,\BSpin 3] &\cong  [X, \BSp1]\cong \rksp(X), \\
    [X,\BSpin 4] &\cong  [X, \BSp1\times \BSp 1]\cong \rksp(X)\oplus \rksp(X), \\
    [X,\BSpin 5] &\cong  [X, \BSp2]\cong \rksp(X), \\
    [X,\BSpin 6] &\cong  [X, \BSU4]\cong [X,\text{BSU}]<\rku(X). 
\end{align*}
}
The invariant $\rho[P]$ can be further simplified by taking the Pontrjagin, Chern, or symplectic Pontrjagin classes, again depending on $n\mod 8$. On the other hand, one can look at the Pontrjagin classes of the bundle $P\times_{\pi_n}\R^n$, where $\pi_n$ is the double cover $\Spin n\to \SO n$. It will be useful for us to compare these classes in dimensions $3\leq n\leq 6$. \\
\lem{\label{lem:cohomo}Let $P$ be a principal $\Spin n$ bundle over base space $X$. \\If $n=3$, then:
$$
    p_1(P\times_{\pi_3}\R^3)=4sp_1 (P\times_{\rho_3}\Ha).
$$
If $n=4$, then, depending on the choice of indices of $\rho_4^1$ and $\rho_4^2$:
\begin{align*}
p_1(P\times_{\pi_4}\R^4)&=2sp_1 (P\times_{\rho_4^1}\Ha) +2sp_1 (P\times_{\rho_4^2}\Ha),\\
    \pm e(P\times_{\pi_4}\R^4)&=sp_1 (P\times_{\rho_4^1}\Ha) -sp_1 (P\times_{\rho_4^2}\Ha).
\end{align*}
If $n=5$, then:

    $$p_1(P\times_{\pi_5}\R^5)=2sp_1 (P\times_{\rho_5}\Ha^2).$$
If $n=6$:
    $$p_1(P\times_{\pi_6}\R^6)=-2c_2 (P\times_{\rho_6^1}\C^4),$$
    $$e(P\times_{\pi_6}\R^6)=c_3 (P\times_{\rho_6^1}\C^4).$$
Here $e$ denotes the Euler class and $sp_i$ denotes the $i$-th symplectic Pontrjagin class, which is the pullback of $(-1)^ic_{2i}$ along the map $\BSp n\to\BU {2n}$ induced by the inclusion $\Sp n\to\U{2n}$.
}
All of the statements above follow from the classical work of Borel and Hirzebruch. We can write down explicit forms of representations $\rho_n$ and $\pi_n$, compute their weights, and use Theorem 10.3 from \cite{BorelHirzebruch} to compare appropriate characteristic classes. The details of these computations can be found in Section \ref{sec:computations}.

\section{\label{sec:main} Classification of real spin bundles over \texorpdfstring{$\CP^3$}{CP3}}
We will need the following:
\lem{\label{lem:stable}(Characteristic classes of stable bundles over $\CP^3$)
\begin{enumerate}
\item The map $\rko(\CP^3)\to H^4(\CP^3;\Z)$ given by the first Pontrjagin class is an isomorphism.
\item The map $\rksp(\CP^3)\to H^4(\CP^3;\Z)$ given by the first symplectic Pontrjagin class is onto and has kernel isomorphic to $\Z/2$.
\item The map $\rku(\CP^3)\to H^*(\CP^3,\Z)$ given by the total Chern class is injective and has in its image every element $1+c_1+c_2+c_3$ with $c_3\equiv0\mod 2$. 
\end{enumerate}
}
\pr{Applying Proposition 2.1 in \cite{Atiyah} we can compute that $\rko(\CP^3)\cong \Z$ and $\rksp(\CP^3)\cong \Z\times \Z/2$. The maps in points (1) and (2) are clearly homomorphisms. Let us denote by $H$ the tautological complex bundle over $\CP^3$. By considering $\SO2$- and $\Sp1\cong \SU2$-bundles of the form $H_\R$ and $H\oplus\overline{H}$ respectively, we can see that the images of the homomorphisms in points (1) and (2) contain a generator of $H^4(\CP^3;\Z)$, hence these homomorphisms are onto. Finally, (3) follows from Theorem A in \cite{66ea2c8e-5263-3d94-93aa-2542242d8d3c}.}
Let $x$ be the generator of $H^2(\CP^3;\Z)$ that is equal to the first Chern class of the tautological bundle.
\thm{\label{thm:main} For each $k,l\in \Z$:
\begin{enumerate}
\item There are exactly two $\Spin 2$ bundles over $\CP^3$ with $p_1=4k^2x^2$.
%\item first
\item There are exactly two $\Spin 3$ bundles over $\CP^3$ with $p_1=4kx^2$.
\item There are exactly four $\Spin 4$ bundles over $\CP^3$ with $p_1=2kx^2, e=lx^2$.
\item There are exactly two $\Spin 5$ bundles over $\CP^3$ with $p_1=2kx^2$.
\item There is exactly one $\Spin 6$ bundle over $\CP^3$ with $p_1=2kx^2,e=2lx^3$.
\item For $n\geq 7$ there is exactly one $\Spin n$ bundle over $\CP^3$ with $p_1=2kx^2$.
\end{enumerate}}
\pr{(1) Since $H^1(\CP^3;\Z_2)=0$,  there is a correspondence between $\Spin 2$-bundles over $\CP^3$ and $\SO 2$-bundles with even Euler class. It follows that for each $k$ there is exactly one $\SO 2$-bundle over $\CP^3$ admitting a spin structure with $e=2k\alpha$, hence there are exactly two $\Spin 2$ bundles over $\CP^3$ with $p_1=e^2=4k^2\alpha^2$. \\
\indent (2) - (5) By Lemma \ref{lem:classification} if $3\leq n\leq 6$, then the invariant $\rho[-]$ defines a correspondence between $\Spin n$-bundles and stable bundles in appropriate $K$-theoretical groups. Lemma \ref{lem:stable} gives us a description of these groups with respect to characteristic classes. Now an application of Lemma \ref{lem:cohomo} gives us the first Pontrjagin classes of corresponding real spin bundles (and Euler classes if $n=4$ or $6$). \\
\indent (6) Since $\BSO n\to\BSO{n+1}$ is a $7$-equivalence for $n\geq 7$, we are in the stable range. Now, by Lemma \ref{lem:stable} for each element of $H^4(\CP^3;\Z)$ there is exactly one $\SO n$ bundle having this element as $p_1$. Since $p_1\equiv w_2^2\mod 2$ and spin bundles have $w_2=0$, the statement follows.

}

\section{\label{sec:app} Construction of a smooth embedding of \texorpdfstring{$\RP^7$}{RP7} in \texorpdfstring{$\RP^{11}$}{RP11}}
We will translate our problem into a problem about bundles over $\CP^3$. For this, we need:
\prop{Let $H$ denote the tautological bundle over $\CP^3$, and let $\eta$ denote the underlying real bundle of $H\otimes_\C H$, then the sphere bundle of $\eta$ is diffeomorphic to $\RP^7$.}
\pr{ Taking out the zero section from $H$ we get the fibration $\C\setminus\{0\}\to\C^4\setminus\{0\}\to\CP^3$. Equip $\C^4$ with an $i$-invariant real inner product. It induces a real inner product on $H$ (and hence on $H\otimes_\C H$), such that the total space of the unit sphere bundle of $H$ is equal to the unit sphere $S^7\subset\C^4$. Now the tensor product $H\to H\otimes_\C H$ sends this $S^7$ onto the total space of the unit sphere bundle of $H\otimes_\C H$. Two points have coinciding images under this map exactly if they are antipodal and therefore the total space of the sphere bundle of $H\otimes_\C H$ is diffeomorphic to the quotient of $S^7$ by antipodism. }

Now we can try to construct an embedding $\CP^3\subset \R^{11}$ such that its normal bundle will have $\eta$ as a summand. In fact\footnote{see Theorem 1.2 in \cite{James}. Another proof relying on the result of Milgram from \cite{Milgram} can be found in \cite{Mukherjee}.}, $\CP^3$ can be embedded smoothly in $\R^9$. Let us denote the normal bundle of this embedding as $N$. Unfortunately, a calculation on characteristic classes shows that the smallest codimension of an Euclidean embedding of $\CP^3$ which might have $\eta$ as a subbundle of its normal bundle is equal to $5$. Nevertheless, we can use inclusion $\R^9\subset\R^{11}$ to obtain an embedding $\CP^3\subset\R^{11}$ having $N\oplus2$ as its normal bundle. Can we now find a bundle $E$ such that $ \eta\oplus E\cong N\oplus2$? Let us, as in Section \ref{sec:main}, use $x$ to denote $c_1(H)$. Standard computations show that:
$$p_1(\eta)=4x^2=p_1(\CP^3)=-p_1(N),$$
$$w_2(\eta)=w_2(\CP^3)=w_2(N)=0.$$
It follows that $E$ has to satisfy $w_2(E)=0$ and $p_1(E)=-8x^2$ (in particular, $E$ admits a spin structure, which is unique, since $H^1(\CP^3;\Z)=0$). From our classification in Theorem \ref{thm:main} we see that there are exactly two bundles $E_0, E_1$ that satisfy these conditions. We will distinguish them using:

\lem{\label{lem:distinguish}Exactly one of the bundles $E_i$ has $\rho[E_i]$ divisible by two.}
\pr{From the proof of Lemma \ref{lem:stable} we know that $\rksp(\CP^3)\cong\Z\oplus \Z_2$ and that the canonical (up to sign) surjection $\Z\times\Z/2\to\Z$ is given by the first symplectic Pontrjagin class. By Lemma \ref{lem:cohomo} we have that $sp_1(\rho[E_i])=\frac{1}{4}p_1(E_i)=2x^2$. Now we see that regardless of the identification of $\rksp(\CP^3)$ with $\Z\times\Z/2$ the elements $\rho[E_i]$ must be of the form $( \pm2, 0)$ and $(\pm2,1)$ and only the first one of them is divisible by two in $\Z\times\Z/2$.}
We will take $E$ to be the bundle $E_i$ having $\rho[E_i]$ divisible by two. We want to show that:

\lem{\label{lem:eq}$E\oplus \eta\cong N\oplus 2$.}
Before we can prove this, we need to make some further observations. Looking at $w_2$ we see that both bundles $E\oplus\eta$ and $N\oplus 2$ admit unique spin structures, so according to Lemma \ref{lem:classification} they are equal if and only if $\rho[E\oplus \eta]=\rho[N\oplus 2]$. In order to compute these invariants, we will need the following:

\lem{\label{lem:tensor} Let $\xi_3$ and $\xi_2$ be respectively $3-$ and $2-$dimensional oriented real vector bundles over some base space $X$ with $H^1(X;\Z_2)=0$. Assume that $\xi_3$ and $\xi_2$ admit (unique) spin structures. Then $\xi_3\oplus\xi_2$ also admits an unique spin structure and, using the conventions stated below, there is a following isomorphism of $\Sp 2$-bundles: $$\pspin{\xi_3\oplus\xi_2}\otimes_{\rho_5}\Ha^2\cong (\pspin{\xi_3}\otimes_{\rho_3}\Ha)\otimes_\R(\pspin{\xi_2}\otimes_{\phi}\R^2).$$
Here $\phi$ denotes the isomorphism $\Spin 2\to\SO 2$, $\pso{-}$ denotes the principal $\SO n$-bundle of ortogonal frames and $\pspin{-}$ denotes the corresponding principal $\Spin n$-bundle. We treat the bundle on the right hand side as a quaternionic $2$-bundle due to the fact that the Kronecker product of matrices gives us a homomorphism $\Sp 1\times \SO 2\to\Sp 2$.}

This lemma can be proven by writing explicitly the lift of inclusion $\SO3\times\SO 2\to\SO 5$ to the corresponding Spin groups:
$$\Sp1\times\SO 2\cong \Spin 3\times \Spin 2\to\Spin 5\cong \Sp 2$$
and observing that it is equivalent to the Kronecker product of matrices. We defer the details to Section \ref{sec:computations}.
\begin{proof}[Proof of Lemma \ref{lem:eq}]
By Lemma \ref{lem:classification} it enough to show that the $\rksp(\CP^3)$-valued invariants $\rho[E\oplus \eta]$ and $\rho[ N\oplus 2]$  are equal. We proceed similarly to the proof of Lemma \ref{lem:distinguish}. The group $\rksp(\CP^3)\cong\Z\times\Z/2$ has a canonical projection onto the $\Z$-summand determined by the symplectic Pontrjagin class $sp_1$. We know from Lemma \ref{lem:cohomo} that for any $\Spin 5$-bundle $\xi$:
$$sp_1(\xi \times_{\rho_5}\Ha^2)=\frac{1}{2}p_1(\xi\times_{\pi_5}\R^5).$$
By the product formula for Pontrjagin classes:
$$p_1(E\oplus \eta)=p_1(N\oplus 2)=-4x^2,$$
hence the projections of $\rho[E\oplus \eta]$ and $\rho[ N\oplus 2]$ to the $\Z$-summand of $\rksp(\CP^3)$ are even and therefore $\rho[E\oplus \eta]=\rho[ N\oplus 2]$ if and only if these classes are both zero or both non-zero in $\rksp(CP^3)\otimes\Z/2$. We will show that the our bundles satisfy the first case, i.e. both classes $\rho[E\oplus \eta]$ and $\rho[ N\oplus 2]$ are divisible by two.\\
\indent We start with $\rho[ N\oplus 2]$. By the definition of the $\rho[-]$-invariant we have that it is equal to the $\rksp$-class of the bundle:
$$\pspin{N\oplus2}\times_{\rho_5}\Ha^2.$$
Using Lemma \ref{lem:tensor} we have following isomorphisms of $\Sp2$-bundles:
\begin{align*}\pspin{N\oplus2}\otimes_{\rho_5}\Ha^2&\cong (\pspin{N}\otimes_{\rho_2}\Ha)\otimes_\R(\pspin{2}\otimes_{\rho_3}\C) \\
& \cong (\pspin{N}\otimes_{\rho_2}\Ha)\otimes_\R 1_{\C}\\
&\cong( \pspin{N}\otimes_{\rho_2}\Ha)\oplus (\pspin{N}\otimes_{\rho_2}\Ha),&
\end{align*}
hence $\rho[N\oplus2]$ is divisible by two in $\rksp(\CP^3)$. \\
\indent Now we treat $\rho[E\oplus \eta]$. Recall that we picked $E$ such that $\rho[E]$ is divisible by two and recall also that by Lemma \ref{lem:classification} quaternionic bundles over $\CP^3$ are stable, i.e. every quaternionic bundle over $\CP^3$ can be uniquely split as a sum of a line bundle and a trivial bundle. We can conclude that there exists bundle $\xi$ satisfying:
$$\xi\oplus \xi=(\pspin{E}\otimes\Ha)\oplus1.$$
Using Lemma \ref{lem:tensor} and additivity of the tensor product we have that in $\rksp(\CP^3)$:
\begin{align*}[\pspin{E\oplus \eta}\otimes_{\rho_5}\Ha^2]&=[ (\pspin{E}\otimes_{\rho_2}\Ha)\otimes_\R(\pspin{\eta}\otimes_{\rho_3}\C)] &\\ &=2[\xi\otimes_\R(\pspin{\eta}\otimes_{\rho_3}\C)]-[1_\Ha\otimes_\R(\pspin{\eta}\otimes_{\rho_3}\C)].\end{align*}
It remains to prove that $[1_\Ha\otimes_\R(\pspin{\eta}\otimes_{\rho_3}\C)]$ is divisible by two in $\rksp(\CP^3)$. Consecutive applications of Lemma \ref{lem:tensor} show that:
\begin{align*}[1_\Ha\otimes_\R(\pspin{\eta}\otimes_{\rho_3}\C)]&=\rho[3\oplus\eta]\\&=\rho[(\eta\oplus1)\oplus 2]\\
&=[(\pspin{\eta\oplus 1}\otimes_{\rho_3}\Ha)\otimes_\R(\pspin2\otimes_{\rho_2}\C)]\\
&= [(\pspin{\eta\oplus 1}\otimes_{\rho_3}\Ha)\otimes_\R1_\C]\\
&=2[(\pspin{\eta\oplus 1}\otimes_{\rho_3}\Ha)].\end{align*}
As a result:
\begin{align*}\rho[E\oplus\eta]=[\pspin{E\oplus \eta}\otimes_{\rho_5}\Ha^2]=2[\xi\otimes_\R(\pspin{\eta}\otimes_{\rho_3}\C)]-2[(\pspin{\eta\oplus 1}\otimes_{\rho_3}\Ha)].\end{align*}
We conclude, that $\rho[E\oplus\eta]$ is also divisible by two, which finishes the proof.
\end{proof}
We now see that $\eta\cong(H\otimes_{\C}H)_\R$ is a subbundle of the bundle $N+2$, which is normal to an embedding of $\CP^3$ in $\R^{11}$. In particular, the sphere bundle of $(H\otimes_\C H)_\R$ embeds smoothly in $\R^{11}$. Thus, we have proven:
\thm{There is a smooth embedding of $\RP^7$ into $\R^{11}$.}

\section{\label{sec:computations} Proofs of lemmas \ref{lem:cohomo} and \ref{lem:tensor}}

 In order to prove lemmas \ref{lem:cohomo} and \ref{lem:tensor}, we have to write explicitly exceptional isomorphisms of spin groups. We will describe them using the action of $\GL{4,\C}$ on $\Lambda^2(\C^4)$. We choose a hermitian metric on $\C^4$. It allows us to define a symplectic form $\omega$ and a hermitian metric $\langle - ,-  \rangle$ on $\Lambda^2(\C^4)$. The subgroup $\SU 4$ preserves both $\langle -,-  \rangle$ and $\omega$, hence it preserves the decomposition of $\Lambda^2(\C^4)$ into $6$-dimensional real subspaces $\Lambda^+, \Lambda^-$ consisting of forms satisfying $\alpha=\ast\overline{\alpha}$ and $\alpha=-\ast\overline{\alpha}$ respectively. These subspaces are equipped with a positive-definite real inner product equal to the real part of $\langle - ,-  \rangle$. We choose to work with $\Lambda^-$ and define its ortonormal basis:
\begin{align*}
&\omega_1=ie_1\wedge e_2 + ie_3\wedge e_4, \\
& \omega_2=e_1\wedge e_2 - e_3\wedge e_4, \\
& \omega_3=ie_1\wedge e_3 - ie_2\wedge e_4, \\
& \omega_4=e_1\wedge e_3 + e_2\wedge e_4, \\
& \omega_5=ie_1\wedge e_4 + ie_2\wedge e_3, \\
& \omega_6=e_1\wedge e_4 - e_2\wedge e_3, 
\end{align*}
For any $I\subseteq {1,\dots,6}$ the subgroup of $\SU 4$ stabilizing $\omega_i:i\in I$ acts on the forms $\omega_i:i\not\in I$. In this way, we can produce an explicit form of the exceptional isomorphisms of the group $\Spin{n-|I|}$. Note that these stabilizers can be computed by writing matrices $M_i$ corresponding to the forms $\omega_i$ in the basis $e_j$ and using the fact that the condition:w
$$U^\intercal M_iU=M_i$$
is equivalent to the following one, which is linear in the entries of $U$:
\begin{equation}
    M_iUM_i^{-1}=\ov{U}.\tag{$\ast$}
\end{equation}

\subsection{\label{subsec:cohomo} Proof of Lemma \ref{lem:cohomo}}
\begin{proof}
By the functoriality of characteristic classes, it is enough to work with classifying spaces. We start by recalling the general theory described in sections 9 and 10 of \cite{BorelHirzebruch}. Suppose that we are given a compact Lie group $G$, its maximal torus $T$ and a homomorphism $\phi$ from $G$ to some group of the form $\SO k, \U k$ or $\Sp k$ . In particular, $\phi$ induces a map $B\phi$ of classifying spaces. Let $pr$ denote the projection in the fiber bundle:
$$G/T\hookrightarrow EG/T\xrightarrow{pr} BG.$$
Theorem 10.3 from \cite{BorelHirzebruch} gives us a way to calculate pullbacks of total Pontrjagin, Chern or symplectic Pontrjagin classes (denoted as $p, c$ and $sp$ respectively) along $B\phi\circ pr$ in terms of the weights of our homomorphism. These weights are elements of $V^*$, the dual space of the universal cover of $T$, and hence can be interpreted as elements of $H^1(T;\Z)$. We can now use the transgression in the spectral sequence associated with the bundle: 
$$T\hookrightarrow EG \to EG/T$$
to identify weights with elements of $H^2(EG/T;\Z)$. Following the authors of \cite{BorelHirzebruch} we will denote weights and corresponding cohomology classes with the same symbols. \\
\indent From now on, we take $G$ to be isomorphic to $\Spin n$ for some $3\leq n\leq6$. Having this assumption, we can use the fact that $G$ is isomorphic to one of the groups $\Sp1, \Sp 1\times \Sp 1, \Sp 2$ or $\SU 4$ to see that the map:
$pr^*:H^*(BG;\Z)\to H^*(EG/T;\Z)$ is injective. As a result, in order to compare pullbacks of characteristic classes to $H^*(BG;\Z)$ it is enough to compare them in $H^*(EG/T;\Z)$. We will now describe our computations in more details. \\
\indent 

\indent ($n=3$) Checking the condition $(\ast)$ shows that the stabilizer of forms $\omega_1, \omega_2, \omega_6$ consists of matrices of the form:
$$\begin{pmatrix}
z & -\ov{w} & 0 & 0 \\
w & \ov{z} & 0 & 0 \\
0 & 0 & z & -\ov{w} \\
0 & 0 & w & \ov{z} \\
\end{pmatrix},$$
where $z,w\in \C, zw\neq0$. It acts on the forms $\omega_3, \omega_4, \omega_5$ and we can use them as a basis for the matrix representation $\Spin 3\xrightarrow{\pi_3}\SO 3$. This representation maps the diagonal maximal torus $T$ to the standard maximal torus of $\SO 3$ and has weight equal to $\pm2x$, where $x$ is dual to a unit vector in the universal cover of $T$. By Theorem 10.3 from \cite{BorelHirzebruch}:
$$pr^*(p(\ESpin 3\times_{\pi_3}\R^3))=1+4x^2.$$
Now in order to compute the total symplectic Pontrjagin class of $\ESpin 3\times_{\rho_3}\Ha$ it seems natural to consider the representation $\Spin 3\xrightarrow{\rho_3}\Sp 1\to \U 2$. We will instead proceed directly using the fact that $\rho_3$ is an isomorphism to apply equation (3) from Section 9.6 of \cite{BorelHirzebruch} and get:
$$pr^*(sp(\ESpin 3\times_{\rho_3} \Ha))=1+x^2.$$
Comparing elements of degree $4$ proves our assertion. \\

\indent (n=4) Stabilizer of the forms $\omega_1$ and $\omega_2$ is isomorphic to $\Spin 4$ and consists of elements of $\SU4$ of the form:
$$\begin{pmatrix}
z_1 & -\ov{w_1} & 0 & 0 \\
w_1 & \ov{z_1} & 0 & 0 \\
0 & 0 & z_2 & -\ov{w_2} \\
0 & 0 & w_2 & \ov{z_2} \\
\end{pmatrix}$$
We can again pick the diagonal maximal torus and equip its universal cover with coordinates $x_1, x_2$ corresponding to maximal tori of upper and lower $2\times 2$ blocks respectively. This torus is mapped via the action on $\omega_i:i=3,\dots,6$ into the standard maximal torus of $\SO4$ and this representation hainals weights $-x_1-x_2$ and $-x_1+x_2$. Again using Theorem 10.3 from \cite{BorelHirzebruch} we obtain:
$$pr^*(p(\ESpin 4\times_{\pi_4}\R^4))=(1+(-x_1-x_2)^2)(1+(-x_1+x_2)^2).$$
$$pr^*(e(\ESpin 4\times_{\pi_4}\R^4))=(-x_1-x_2)(-x_1+x_2)=x_1^2-x_2^2.$$
On the other hand the representations $\rho_4^i:\Spin 4\to \Sp 1$ for $i=1,2$ have weights equal to $x_1$ and $x_2$ and hence by Section 9.6 from \cite{BorelHirzebruch}:
$$pr^*(sp(\ESpin 4\times_{\rho_4^i}\Ha))=1+x_i^2.$$
Comparing terms of degree $4$ we get our result. Note that we could switch the order of $x_1$ and $x_2$, hence the ambiguity regarding the sign of the Euler class. \\

\indent (n=5) Stabilizer of the form $\omega_1$ is isomorphic to $\Spin 5$ and consists of elements of $\SU4$ of the form:
$$\begin{pmatrix}
z_1 & -\ov{w_1} & z_2 & -\ov{w_2} \\
w_1 & \ov{z_1} & w_2 & \ov{z_2} \\
z_3 & -\ov{w_3} & z_4 & -\ov{w_4} \\
w_3 & \ov{z_3} & w_4 & \ov{z_4} \\
\end{pmatrix}$$
We pick the same maximal torus as in the case $n=4$, equip its universal cover with the same coordinates and get the same weights of representation $\pi_5$. We have again that:
$$pr^*(p(\ESpin 5\times_{\pi_5}\R^5))=(1+(-x_1-x_2)^2)(1+(-x_1+x_2)^2).$$
This time, however Section 9.6 from \cite{BorelHirzebruch} tells us that:
$$pr^*(sp(\ESpin 5\times_{\rho_5}\Ha^2))=(1+x_1^2)(1+x_2^2).$$
We again compare terms of degree $4$ and get our result. \\

\indent (n=6)
We pick a diagonal maximal torus in $\Spin 6\cong \SU4$ and equip its universal cover with coordinates $x_1, x_2, x_3$ such that the projection of $(x_1,x_2,x_3)$ to $T^3$ has the form:
$$\begin{pmatrix}
e^{2\pi(x_1+x_2+x_3)} & 0 & 0 & 0 \\
0 & e^{-2\pi x_1} & 0 & 0 \\
0 & 0 & e^{-2\pi x_2} & 0 \\
0 & 0 & 0 & e^{-2\pi x_3} \\
\end{pmatrix}.$$
Now a direct computations shows that the weights of $\pi_6$ are: $(-x_2-x_3), (-x_1-x_3),$ and $(-x_2-x_3).$ An application of Theorem 10.3 from \cite{BorelHirzebruch} shows that:
$$pr^*(p(\ESpin 6\times_{\pi_6}\R^6))=(1+(-x_2-x_3)^2)(1+(-x_1-x_3)^2)(1+(-x_2-x_3)^2),$$
$$pr^*(e(\ESpin 6\times_{\pi_6}\R^6))=(-x_2-x_3)(-x_1-x_3)(-x_2-x_3).$$
We can also apply Theorem 10.3 to the inclusion $\Spin 6\xrightarrow{\cong} \SU 4\xhookrightarrow{i}\U 4$ and get:
$$pr^*(c(\ESpin 6\times_{\rho_6}\C^4))=(1+x+y+z)(1+x)(1+y)(1+z).$$
Comparing terms of degrees $4$ and $6$ we get our equalities.
\end{proof}

\subsection{Proof of Lemma \ref{lem:tensor}}
\begin{proof}
Consider the following diagram of groups:
\begin{center}
\begin{tikzcd}
\Sp{1}\times\SO{2}  \arrow[dashrightarrow]{r}{\overline{f}}& \Sp{2} \\
\Spin{3}\times\Spin{2}\arrow[u, "\rho_{3}\times\rho_{2}"', "\cong"] \arrow{d}{\pi_{3}\times\pi_{2}} \arrow[dashrightarrow]{r}{f} & \Spin{5} \arrow[u, "\rho_5"', "\cong"]  \arrow{d}{\pi_{5}} \\
\SO{2}\times\SO{3} \arrow[hook]{r}{i} & \SO{5}.
\end{tikzcd}
\end{center}
Since $\Spin3\times\Spin 2$ is $1$-connected, there are unique arrows $f$ and $\overline{f}$ that make this diagram commute. We can place the bundles $\xi$ and $\eta$ in a corresponding diagram for classifying spaces: 
\begin{center}
\begin{tikzcd}
& & & \BSp{1} \times \BSO{2}\arrow[dashrightarrow]{r}{B\overline{f}}  & \BSp{2} & & &\\
&  & &\BSpin{3}\times\BSpin{2} \arrow[u, "B\rho_{3}\times B\rho_{2}"', "\cong"] \arrow{d}{B\pi_{2}\times B\pi_{3}} \arrow[dashrightarrow]{r}{Bf} & \BSpin{5} \arrow[u, "B\rho_5"', "\cong"]\arrow{d}{B\pi_{5}} & & & \\
X\arrow[urrr, "\pspin{\xi}\times\pspin{\eta}" near end] \arrow[rrr, "\pso{\xi}\times\pso{\eta}"']  & & &\BSO{3}\times\BSO{2} \arrow[hook]{r}{Bi} & \BSO{5}.& & &
\end{tikzcd}
\end{center}
Since $\pso{\xi\oplus\eta}\sim i\circ(\pso\xi\times\pso\eta)$ admits a unique lift to $\BSpin 5$, we have:
\begin{align*} \pspin{\xi\oplus\eta}&\sim Bf\circ(\pspin\xi\times\pspin\eta),\\
B\rho_5\circ\pspin{\xi\oplus\eta}&\sim B\overline{f}\circ (B\rho_3\times B\rho_2)\circ(\pspin\xi\times\pspin\eta).\end{align*}
Therefore, in order to compare $\pspin{\xi\otimes\eta}\times_{\rho_5}\Ha^2$ with $\pspin{\xi}\times_{\rho_3}\Ha$ and $\pspin{\eta}\times_{\rho_2}\C$, it suffices to understand $B\overline{f}$. \\
\indent We can think of $\Spin 5\cong \Sp 2$ as a subgroup of $\Spin 6\cong\SU 4$ stabilizing $\omega_1$. We have already seen in Subsection \ref{subsec:cohomo} that it consists of matrices of the form:
$$\begin{pmatrix}
z_1 & -\ov{w_1} & z_2 & -\ov{w_2} \\
w_1 & \ov{z_1} & w_2 & \ov{z_2} \\
z_3 & -\ov{w_3} & z_4 & -\ov{w_4} \\
w_3 & \ov{z_3} & w_4 & \ov{z_4} \\
\end{pmatrix}, \quad z_i, \ w_i\in\C\text{ for }i=1,2,3,4.$$
We can also check the condition $(\ast)$ to see that the subgroups $\Spin 3\cong \Sp 1$ and $\Spin 2\cong \SO 2$ stablilizing $\omega_1, \omega_2, \omega_6$ and $\omega_1, \omega_3,\omega_4, \omega_5$ respectively are of the following form:
$$\begin{pmatrix}
z & -\ov{w} & 0 & 0 \\
w & \ov{z} & 0 & 0 \\
0 & 0 & z & -\ov{w} \\
0 & 0 & w & \ov{z} \\
\end{pmatrix}, \ z,w\in\C, \quad
\begin{pmatrix}
a & 0 & -b & 0 \\
0 & a & 0 & -b \\
b & 0 & a & 0 \\
0 & b & 0 & a \\
\end{pmatrix},\ a,b\in\R.
$$
It follows that the map $\overline{f}$ is just the complex Kronecker product of matrices, which is the same as the real Kronecker product, since the entries of matrices from $\SO 2$ in the form above are real. This proves our statement.
\end{proof} 

\bibliographystyle{alpha}
\nocite{*}
\bibliography{citations_1}
\end{document}